\documentclass[reqno]{amsart}
\usepackage[margin=3.5cm]{geometry}
\usepackage{xcolor}
\usepackage{amsfonts,amsthm,amssymb,verbatim,amscd,amsmath,blindtext}
\usepackage{multirow}
\usepackage{graphicx}
\usepackage{enumitem}
\usepackage{amssymb}
\usepackage{ytableau}
\usepackage{tikz}
\usepackage{tikz-cd}
\usepackage{float,pifont}
\usepackage{mathtools}
\usepackage[pagebackref,colorlinks=true,citecolor=blue]{hyperref}
\usepackage[capitalise]{cleveref}

\newtheorem{theorem}{Theorem}[section]

\newtheorem{lemma}[theorem]{Lemma}

\newtheorem{proposition}[theorem]{Proposition}

\theoremstyle{definition}

\DeclareMathOperator{\lcm}{lcm}
\DeclareMathOperator{\mingens}{Mingens}

\DeclareMathOperator{\pd}{pd}
\DeclareMathOperator{\reg}{reg}

\begin{document}

\title[Homological invariants of polarized neural ideals]{From neural codes to homological invariants: regularity and projective dimension of\\
polarized neural ideals}

\author[T. Chau]{Trung Chau}
\address{Chennai Mathematical Institute, Siruseri, Tamil Nadu 603103. India}
\email{chauchitrung1996@gmail.com}


\subjclass[2020]{}

\begin{abstract}
   Neural codes form an algebraic framework to study the nervous system, and understanding  neural codes is a key goal
 of mathematical neuroscience. Neural rings and ideals are the tools connecting neuroscience and commutative algebra. In this article, we study the projective dimension and (Castelnuovo-Mumford) regularity of polarized neural ideals on $n$ neurons. Particularly, we find all the possible values for these two invariants. Moreover, we characterize when these ideals have linear resolution or linear quotients, assuming that they are generated in degree $n$.
\end{abstract}

\maketitle

\section{Introduction}

In 1943, McCulloch and Pitts proposed a mathematical model to study the nervous system by focusing on the ``all-or-none" nature of neuron-firing. Since then, the combinatorial, topological, geometric, and coding-theoretical structure of neural responses have been examined from a wide range of perspectives \cite{neural-4,neuron-2,neuron-3,neuro-0,neural-ideals-canonical,neural-5,neuron-1}. We will focus on the algebraic side. In 2013,  Curto et. al. introduced the concept of neural rings, bridging the fields of neuroscience and coding theory via algebraic tools and methods. 

Given $n$ neurons. Each neuron either fires or not, and this ``all-or-none" nature allows researchers to associate these $n$ neurons with a binary codeword $(c_1,\dots, c_n)\mathcal{C}\subseteq \{0,1\}^n$ of length $n$. Here $c_i=0$ means the $i$-th neuron is not firing. Now one can define the associated \emph{neural ring} to be $\mathbb{F}_2[x_1,\dots, x_n]/I_\mathcal{C}$ where $I_\mathcal{C}$ is the vanishing ideal of $\mathcal{C}$. The authors of \cite{neuro-0} then defined a \emph{neural ideal} $J_\mathcal{C}\subset I_\mathcal{C}$ such that
\[
I_\mathcal{C}=J_\mathcal{C}+ (x_i(1-x_i)\mid i\in [1,n]).
\]
The idea is that any relations $x_i(1-x_i)$ are trivial, as it holds for any neural code, and by eliminating them from the set of generators, one obtains a smaller ideal to study, namely $J_\mathcal{C}$. This neural ideal faithfully captures the combinatorial constraints inherent to the code while providing an interface with commutative algebra and algebraic geometry. The perspective proved fruitful, linking neural coding theory to primary decomposition of neural ideals and  questions regarding convexity, and establishing a bridge between the structure of neural codes and the algebraic invariants of the ideals they generate.

Neural ideals on $n$ neurons are generated by \emph{pseudomonomials}, i.e., polynomials of the form $\prod_{i\in \sigma} x_i \prod_{j\in \tau}(1-x_j)$ for some subsets $\sigma$ and $\tau$ of $\{1,\dots, n\}$. Most of these polynomials are not homogeneous under any grading (excluding the grading where the variables are of degree 0). This poses a challenge for researchers to directly apply the many classical tools in algebra and geometry. To remedy this, G\"unt\"urk\"un et. al. \cite{neural-5} introduced the concept of polarization to transform a neural ideal into a \emph{polarized neural ideal}, a squarefree monomial ideal in $\mathbb{F}_2[x_1,\dots, x_n,y_1,\dots, y_n]$. This opens the door for the many framework that can be imposed to study neural ideals, e.g., Stanley-Reisner theory. Roughly speaking, each pseudomonomial generator $\prod_{i\in \sigma} x_i \prod_{j\in \tau}(1-x_j)$ is ``polarized" into $\prod_{i\in \sigma} x_i \prod_{j\in \tau}y_j$, a squarefree monomial, and the ``polarized" monomials generate what we call \emph{polarized neural ideals}. These ideals thus serve as a gateway to the rich combinatorial and homological machinery of squarefree monomial ideals. For example, the minimal primes of a polarized neural ideal correspond to facets of a simplicial complex encoding the code’s combinatorics, and its resolutions can be studied using cellular and combinatorial methods unavailable in the original Boolean setting. This development reflects a broader trend in the algebraic study of combinatorial data: translating problems into squarefree settings to leverage the geometry and topology of simplicial complexes and the power of homological invariants. We refer to \cite{neural-5} for more.

Surprisingly, there have not been much progress on this algebraic side. In this article, we continue the study of  polarized neural ideals from an algebraic perspective. Each (monomial) ideal has a minimal free resolution, which encodes all the algebraic and geometric information of the ideal. In general, an explicit construction of the minimal free resolution is a difficult task, and is often accomplished for only special classes of ideals (even for monomial ideals) \cite{BM20, BW02, BPS98, BS98, CK24, CT2016, CEFMMSS21, CEFMMSS22, Ly88, OY2015, Vel08}. However, to study related invariants, there are more tools to offer. The invariants we will work on in this article are \emph{projective dimension} and \emph{(Castelnuovo-Mumford) regularity} (we refer to Section~\ref{sec:prem} for unexplained terminology). Roughly speaking, these two invariants measure the ``size" of a minimal free resolution. For more history on these invariants, we refer to some recent work and survey \cite{DHS,regularity-survey}.

In this article, by a \emph{polarized neural ideal on $n$ neurons}, we mean a squarefree monomial ideal in $\Bbbk[x_1,\dots, x_n,y_1,\dots y_n]$ such that $x_iy_i$ does not divide any of its minimal monomial generator for any $i\in [1,n]$, where $\Bbbk$ is any field and $n$ is any positive integer. This is more general than the polarized neural ideals in literature, but these conditions are common in combinatorial commutative algebra. Note that the minimal monomial generators of a polarized neural ideal can only be at most $n$. Let $\pd I$ and $\reg I$ denote the projective dimension and regularity of $I$, respectively, for a polarized neural ideal $I$.  Our first result is on the possible values of these invariants of $I$, where $I$ is a neural ideal on $n$ neurons (with or without some conditions). Effectively, this gives sharp lower and upper bounds on these invariants.

\begin{theorem}[\protect{Theorems~\ref{thm:naive-bounds} and \ref{thm:naive-2}}]
    Let $I$ be a nonzero polarized neural ideal on $n$ neurons. Then $\{0,\dots, 2n-1\}$ and $\{1,\dots, 2n-1\}$ are the sets of all possible values of $\pd I$ and $\reg I$, respectively. If $I$ is additionally assumed to be generated in degree $n$, the sets become $\{0,\dots, n\}$ and $\{n,\dots, 2n-1\}$, respectively. 
\end{theorem}

It is a question of interest of when a squarefree monomial ideal has linear resolution or linear quotients \cite{AJM2024,BJT19,Fai17,Fic25, HS25,JahanZheng2010,SharifanVarbaro2008} (we again refer to Section~\ref{sec:prem} for unexplained terminology). The motivation for linear resolution is the classical result by Fr\"oberg, which states that the edge ideal of a finite simple graph has linear resolution if and only if the graph is co-chordal, bridging the fields of graph theory and commutative algebra. On the other hand, the property of having linear quotients is intricately linked to the shellability of certain simplicial complexes. A necessary condition for a monomial ideal to have linear resolution is that all monomial generators are generated in the same degree. In the case where polarized neural ideals on $n$ neurons are generated in the maximum degree $n$, we obtain a recursive classification based on the properties of having linear resolution and linear quotients.

\begin{theorem}[\protect{Theorems~\ref{thm:linear-resolution/quotients-inductive} and \ref{thm:LR-LQ-the-same}}]
    Let $I$ be a polarized neural ideal on $n$ neurons and assume that $I$ is generated in degree $n$. Then $I = x_nJ+y_nK$ where $J,K$ are (uniquely determined) neural ideals in $n-1$ neurons, and in particular do not involve $x_n$ or $y_n$. Then the following are equivalent:
    \begin{enumerate}
        \item $I$ has linear resolution;
        \item $I$ has linear quotients;
        \item $J$ and $K$ have linear resolution, and one contains the other;
        \item $J$ and $K$ have linear quotients, and one contains the other.
    \end{enumerate}
\end{theorem}

By placing polarized neural ideals in the broader history of combinatorial commutative algebra and neural coding theory, this work advances both the algebraic understanding of neural codes and the catalog of monomial ideals whose homological invariants are controlled by combinatorial data. It opens the door to further connections between neuroscience-inspired combinatorics and classical commutative algebra.

The article is structured as follows. We provide the commutative background in Section~\ref{sec:prem}. In Section~\ref{sec:bounds}, we obtain the bounds as well as all the possible values for projective dimension and regularity of neural ideals. Finally, we study the properties of having linear resolution and linear quotients in Section~\ref{sec:linear-res}.

\subsection*{Acknowledgments} The author is supported by the Infosys Foundation. The author would like to thank Hugh Geller and Rebecca R.G. for inspiring talks leading to this article. This article was an idea in Spring 2024, and the author would like to dedicate it to Deborah Anne Wooton, for everything they did.

\section{Preliminaries}\label{sec:prem}

In this article, for a positive integer $n$, we use $[n]$ to denote the set $\{1,\dots, n\}$. Throughout the article, $\Bbbk$ denotes a field.
	
	\subsection{Linear quotients and linear resolution}
	
	Let $S$ be the polynomial ring $\Bbbk[x_1,\dots, x_n]$ over a field $\Bbbk$, $\mathfrak{m}$ its irrelevant ideal, and $M$ a finitely generated module over $S$. A \emph{free resolution} of $M$ over $S$ is a complex of free $S$-modules
	\[
	\mathcal{F}\colon \cdots \to F_r\xrightarrow{\partial} F_{r-1} \to \cdots \to F_1 \xrightarrow{\partial} F_0\to 0
	\]
	such that $H_0(\mathcal{F})\cong M$ and $H_i(\mathcal{F})=0$ for any $i>0$. Moreover, $\mathcal{F}$ is called \emph{minimal} if $\partial(F_{i+1})\subseteq \mathfrak{m}F_{i}$ for any $i$. It is known that the minimal free resolution of $M$ is finite, i.e., $F_i=0$ if $i\gg 0$. We can thus define the \emph{projective dimension} of $M$, denoted by $\pd M$, to be
    \[
    \pd M \coloneqq \max \{i\mid F_i\neq 0\}.
    \]
	
	In the case where $M$ is $\mathbb{N}$-graded, it is well-known that it has an $\mathbb{N}$-graded minimal free resolution. In other words, the free modules $F_i$ can be given a shift so that the differentials are homogeneous. We can thus set $F_i=\oplus_{j\in \mathbb{N}} S(-j)^{\beta_{i,j}(M)}$ for any integer $i$. The numbers $\beta_{i,j}(M)$ are called \emph{Betti numbers} of $M$. The \emph{(Castelnuovo-Mumford) regularity} of $M$, denoted by $\reg M$, is defined to be
	\[
	\reg M\coloneqq \max \{ j - i \mid \beta_{i,j}(M)\neq 0 \}.
	\]
	In this article we will only consider the case where $M=I$ is a monomial ideal of $S$. Recall that $I$ has a unique minimal monomial generating set, which we denote by $\mingens(I)$. The following is straightforward, but we record this below for easy reference.

    \begin{lemma}\label{lem:reg-max-degree}
        Let $I$ be a monomial ideal. Then $\reg I\geq \max_{m\in \mingens(I)} \{\deg(m)\}$.\qed
    \end{lemma}
    
    If $I$ is generated by monomials in the same degree $d$, We say that $I$ is \emph{generated in degree $d$}. A monomial ideal $I$ generated in degree $d$ is said to have \emph{linear resolution} if $\reg I = d$. We say that $I$ has \emph{linear quotients} if after a relabelling, we have $\mingens(I)=\{ m_1, m_2,\dots, m_q\}$ where the colon ideal
	\[
	(m_1,m_2,\dots, m_k)\colon (m_{k+1}) 
	\]
	is generated by a set of variables of $S$, for each $k\in [q-1]$. We record the formulas for colon ideals (and intersection) below.

    \begin{lemma}[\protect{\cite[Propositions 1.2.1 and 1.2.2]{HerzogHibiBook}}]\label{lem:formulas}
        Let $I,J$ be monomial ideals and $u$ a monomial. Then
        \begin{align*}
            I\cap J&=\left( \lcm(m,m')\mid m\in \mingens(I) \text{ and } m'\in \mingens(J) \right),\\
            I\colon u&= \left( \frac{m}{\gcd(u,m)} \mid m\in \mingens(I)  \right).
        \end{align*}
    \end{lemma}
    
    The following result is well-known.
	
	\begin{lemma}[{\cite[Theorem 8.2.15]{HerzogHibiBook}}]\label{lem:LQ-implies-LR}
	   Monomial ideals generated in a fixed degree with linear quotients have linear resolution.
	\end{lemma}
	
	For a monomial ideal $I$ and a monomial $m$, let $I^{\leq m}$ be the monomial ideal generated by monomials in $\mingens(I)$ that divide $m$. The following is a direct corollary of the well-known Restriction Lemma.
	
	\begin{lemma}[{\cite[Lemma 4.1]{HHZ2004}}]\label{lem:linear-resolution-restriction}
		If a monomial ideal $I$ has linear resolution, so does $I^{\leq m}$ for any monomial $m$.
	\end{lemma}
	
	In general, an analog of the Restriction Lemma for the property of having linear quotients does not hold. However, it does when the ideal is equigenerated.
	
	\begin{lemma}[{\cite[Proposition 2.6]{HMRG20}}]\label{lem:linear-quotients-restriction}
		If an equigenerated monomial ideal $I$ has linear quotients, then so does $I^{\leq m}$ for any monomial $m$.
	\end{lemma}

    We will also recall how the algebraic invariants and properties of $I$ and $mI$ are related, where $I$ is a monomial ideal and $m$ a monomial. The following is likely known to experts.

    \begin{lemma}\label{lem:mI-I}
        Let $I$ be a monomial ideal and $m$ a monomial. Then
        \begin{enumerate}
            \item $I$ has linear resolution or linear quotients if and only if so does $mI$;
            \item $\pd(mI)=\pd I$;
            \item $\reg(mI)=\reg I$.
        \end{enumerate}
    \end{lemma}

    \begin{proof}
        As polynomial rings are domains, we have an isomorphism of modules $mI\cong I(-\deg(m))$. Statements (2) and (3), and (1) in the case of linear resolution, then follow. It remains to show that $I$ has linear quotients if and only if so does $mI$. Set $\mingens(I)=\{m_1,\dots, m_q\}$. Then $\mingens(mI)=\{mm_1,\dots, mm_q\}$ The result then follows from the identity
        \[
        (m_1,\dots, m_i)\colon m_{i+1} = (mm_1,\dots, mm_i)\colon mm_{i+1} 
        \]
        for any $i\in [q-1]$ (Lemma~\ref{lem:formulas}).
    \end{proof}
	
\subsection{Betti splittings}

Let $I,J,K$ be monomial ideals such that $I=J+K$. We say that $I=J+K$ is a \emph{Betti splitting} if 
	\[
	\beta_{i,j}(I)=\beta_{i,j}(J)+\beta_{i,j}(K)+\beta_{i-1,j}(J\cap K) \;\;\;\text{ for any integers } i,j.
	\]

A popular method to find Betti splittings is the following.

\begin{lemma}[\protect{\cite[Corollaries 2.2 and 2.7]{splitting}}]\label{lem:x-splitting}
    Let $I$ be a monomial ideal and $x$ a variable. Then $I=xJ+K$ where $J$ and $K$ are (uniquely determined) ideals generated by monomials not divisible by $x$. If $J$ has linear resolution, then $I=xJ+K$ is a Betti splitting. In particular, we have
    \begin{align*}
        \pd I &= \max \{\pd (xJ), \ \pd K, \ \pd(xJ\cap K) +1 \},\\
        \reg I &= \max \{\reg (xJ), \ \reg K, \ \reg(xJ\cap K) -1  \}.
    \end{align*}
\end{lemma}

Such a Betti splitting is also known as \emph{$x$-splitting}.

Polarized neural ideals are monomial ideals of special forms. Thus we will only use the following special case of Lemma~\ref{lem:x-splitting}

\begin{lemma}\label{lem:x-splitting-special}
    Let $I$ be a monomial ideal and $x,y$ variables. Assume that $I=xJ+yK$ where $J$ and $K$ are generated by monomials not divisible by $x$ or $y$. If $J$ has linear resolution, then $I=xJ+yK$ is a Betti splitting. In particular, we have
    \begin{align*}
        \pd I &= \max \{\pd (J), \ \pd K, \ \pd(J\cap K) +1 \},\\
        \reg I &= \max \{\reg J+1, \ \reg K+1, \ \reg(J\cap K) +1  \}.
    \end{align*}
\end{lemma}

\begin{proof}
    This follows from Lemmas~\ref{lem:mI-I} and \ref{lem:x-splitting}.
\end{proof}

\section{Some bounds on algebraic invariants of polarized neural ideals}\label{sec:bounds}

In this section we will give sharp bounds on the projective dimension and regularity of a polarized neural ideal on $n$ neurons, as well as provide all of their possible values. We start with computing these invariants for some special polarized neural ideals.

Recall from \cite{dominant} that a set $\mathcal{G}=\{m_1,\dots, m_q\}$ of squarefree monomials  is called \emph{dominant} if for each $k\in [q]$, there exists a variable $x_k$ such that $m_k$ is the only monomial in $\mathcal{G}$ that is divisible by $x_k$. We recall the following result.

\begin{lemma}[\protect{\cite[Theorem 4.8 and Corollary 4.9]{dominant}}]\label{lem:reg-pd-dominant}
    Let $I$ be a squarefree monomial ideal such that $\mingens(I)=\{m_1,\dots, m_q\}$ is a dominant set. Then $\pd I=q-1$ and $\reg I = \deg (\lcm(m_1,\dots, m_q))-q+1$.
\end{lemma}

This will be the only tool we need to compute the invariants for the following polarized neural ideals. The proofs are a straightforward application of Lemma~\ref{lem:reg-pd-dominant}. We leave the details to interested readers. 

\begin{proposition}\label{prop:pd-n-neurons}
    Let $n$ be a positive integer. For each $k\in [n]$, set $m_k= x_ky_1y_2\cdots \widehat{y_k}\cdots y_n$. Then for each $k\in [n]$, the ideal $(m_1,\dots, m_k)$ has projective dimension $k-1$. In particular, the projective dimension of a polarized neural ideal on $n$ neurons generated in degree $n$ can achieve any value in $[0,n-1]$. \qed
\end{proposition}

\begin{proposition}\label{prop:reg-n-neurons}
    Let $n$ be a positive integer. For each $k\in [n]$, the ideal generated by two monomials
    \[
    \left(\prod_{i=1}^n x_i\right)\quad \text{ and } \quad  \left(\prod_{i=1}^k y_i\right)\left(\prod_{j=k+1}^nx_j\right)
    \]
    has regularity $n+k-1$.  In particular, the regularity of a polarized neural ideal on $n$ neurons generated in degree $n$ can achieve any value in $[n,2n-1]$. \qed
\end{proposition}

\begin{proposition}\label{prop:all-reg-pd-neural-ideals}
    Let $n$ be a positive integer. For each $i\in[0,2n-1]$ and each $j\in [1,2n-1]$, let $I$ be the ideal generated by a set of $i+1$ variables and $J$ be the principal ideal generated by a monomial of degree $j$. Then $\pd I = i$ and $\reg J=j$. In particular, the projective dimension of a polarized neural ideal on $n$ neurons can achieve any value in $[0,2n-1]$, and the same holds for regularity in $[1,2n-1]$.\qed 
\end{proposition}

Fix a positive integer $n$. Since polarized neural ideals allow variables, their projective dimension and regularity can achieve a few ``naive" lower and upper bounds, as follows.

\begin{theorem}\label{thm:naive-bounds}
    Let $I$ be a nonzero polarized neural ideal on $n$ neurons. Then $\pd I \in [0,2n-1] $ and $ \reg I\in [1,2n-1]$. Moreover, for each integers $j\in [0,2n-1]$ and $k\in [1,2n-1]$, there exist nonzero polarized neural ideals $J$ and $K$ on $n$ neurons such that $\pd J =j$ and $\reg K = k$.
\end{theorem}

\begin{proof}
    Recall that neural ideals are ideals in a polynomial ring $S$ in $2n$ variables. By Auslander-Buchbaum-Serre's theorem, $\pd(S/I)\leq 2n$, or equivalently, $\pd(I)\leq 2n-1$. Regarding regularity, by definition, $\reg(S/I)\geq 0$ (by looking at the first homological degree of its minimal resolution, for example), and thus $\reg I=\reg S/I+1 \geq 1$. For the upper bound,   
    the following is well-known to experts (see, e.g., \cite[Proposition~6.2]{domination}):
    \begin{equation}\label{reg}
        \reg I \leq \max\{ \deg(\lcm(A)) - |A| \mid  \emptyset \neq A\subseteq \mingens(I)\} + 1.
    \end{equation}
    Since $I$ is a squarefree monomial ideal in $2n$ variables, we have $\deg(\lcm(A))\leq 2n$. Therefore, we have
    \[
    \reg I\leq 2n - 1 + 1=2n.
    \]
    It remains to show that $\reg I\neq 2n$. Indeed, suppose that $\reg(I) =2n$. Then (\ref{reg}) implies that there exists a subset $A$ of $\mingens(I)$ such that $\deg(\lcm(A))=2n$ and $|A|=1$. Since $\mingens(I)$ is a set of squarefree monomials in $2n$ variables, we must have $A=\{\prod_{i=1}^n(x_iy_i)\}$. The only squarefree monomial ideal with $\prod_{i=1}^n(x_iy_i)$ as a minimal monomial generator is the principal ideal generated by it. This contradicts the hypothesis that $I$ is a neural ideal. Therefore, $\reg I \leq 2n-1$, as desired. The second statement follows from Proposition~\ref{prop:all-reg-pd-neural-ideals}.
\end{proof}

We also give bounds when the polarized neural ideal on $n$ neurons is generated in degree $n$.

\begin{theorem}\label{thm:naive-2}
    Let $I$ be a nonzero polarized neural ideal on $n$ neurons, and assume that $I$ is generated in degree $n$. Then $\pd I \in [0,n] $ and $\reg I\in [n,2n-1]$. Moreover, for each integers $j\in [0,n]$ and $k\in [n,2n-1]$, there exist nonzero polarized neural ideals $J$ and $K$ on $n$ neurons such that $J,K$ are generated in degree $n$, and that $\pd J =j$ and $\reg K = k$.
\end{theorem}

\begin{proof}
    Since $I$ is a squarefree monomial ideal generated in degree $n$, by \cite[Theorem~3.1]{Hilbert-Alessandroni}, we have
    \[
    \pd I \leq 2n- (n+1)+1=n,
    \]
    as desired. For $\reg I$, the lower bound follows from Lemma~\ref{lem:reg-max-degree}, and the upper bound from Lemma~\ref{thm:naive-bounds}. The second statement would follow from Propositions~\ref{prop:pd-n-neurons} and \ref{prop:reg-n-neurons} if we can point out an example where $\pd I =n$. Indeed, for each $k\in [n]$, set \[
    I_k=(x_1,y_1)(x_2,y_2)\cdots (x_k,y_k).\]
    It suffices to show that $\pd I_n = n$. 
    
    We will in fact prove something stronger: $\pd I_k =k$ for any $k\in [n]$. We proceed by induction on $k$. If $k=1$, then $\pd I_1=\pd (x_1,y_1)=1$ by Lemma~\ref{lem:reg-pd-dominant}. By induction, we can assume that $\pd I_{k-1}=k-1$ for some $k\geq 2$. Observe that $I_k=x_kI_{k-1}+y_kI_{k-1}$ and $I_{k-1}$ has linear resolution (\cite[Theorem~3.1]{regularity-product}). Thus by Lemma~\ref{lem:x-splitting-special}, we have
    \[
    \pd I_k = \max\{\pd I_{k-1},\  \pd I_{k-1},\ \pd(I_{k-1}\cap I_{k-1}) +1 \} = k,
    \]
    as desired. This concludes the proof.
\end{proof}

\section{Polarized neural ideals with linear resolution or linear quotients}\label{sec:linear-res}

In this section we study when a polarized neural ideal $I$ on $n$ neurons has linear resolution or linear quotients, assuming that $I$ is generated in degree $n$. Recall that $I$ is a squarefree monomial ideal in $\Bbbk[x_1,\dots, x_n,y_1,\dots, y_n]$ such that no monomial in $\mingens(I)$ is divisible by $x_iy_i$, for any $i\in [n]$. Couple this with the condition that $I$ is generated in degree $n$, we deduce that each monomial in $\mingens(I)$ is a squarefree monomial that is divisible by exactly one variable in $\{x_i,y_i\}$, for each $i\in [n]$. In particular, each minimal monomial generator of $I$ is divisible by either $x_n$ or $y_n$, but not both. Therefore, there is a unique decomposition $I=x_nJ+y_nK$ where $J$ and $K$ are squarefree monomial ideals generated by monomials not involving $x_n$ and $y_n$. In particular, $J$ and $K$ are neural ideals on $n-1$ neurons, and are generated in degree $n-1$. This is the key to our results in this section.

\begin{theorem}\label{thm:linear-resolution/quotients-inductive}
    Let $I$ be a polarized neural ideal on $n$ neurons and assume that $I$ is generated in degree $n$. Set $I = x_nJ+y_nK$ where $J,K$ are (uniquely determined) ideals generated by monomials that do not involve $x_n$ and $y_n$. Then the following are equivalent:
    \begin{enumerate}
        \item[(i)] $I$ has linear resolution (resp, linear quotients);
        \item[(ii)] $J$ and $K$ have linear resolution (resp, linear quotients), and one contains the other.
    \end{enumerate}
\end{theorem}
    
\begin{proof}
    We first note that each of $J$ and $K$ is either the zero ideal or a polarized neural ideal on $n-1$ neurons and is generated in degree $n-1$ in the latter case.

    \emph{(i) $\implies $ (ii):} Assume that $I$ has linear resolution (resp, linear quotients). It is straightforward that 
    \[
    I^{\leq x_n \prod_{i=1}^{n-1}(x_iy_i)} = x_nJ \quad \text{and}\quad  I^{\leq y_n \prod_{i=1}^{n-1}(x_iy_i)} = y_nK. 
    \]
    Thus by Lemma~\ref{lem:linear-resolution-restriction} (resp, Lemma~\ref{lem:linear-quotients-restriction}), both $J$ and $K$ have linear resolution (resp, linear quotients). 

    It remains to show that $J$ either contain or is contained in $K$. Remark that $I,J,K$ has linear resolution (as linear quotients implies linear resolution by Lemma~\ref{lem:LQ-implies-LR}). Since $J$ has linear resolution, the decomposition $I=x_nJ+y_nK$ is an $x_n$-splitting by Lemma~\ref{lem:x-splitting-special}. Therefore
    \begin{equation}\label{reg-equa}
        n=\reg I = \max\{ \reg J+1 , \reg K+1, \reg (J\cap K) +1 \} = \max\{ n,\reg(J\cap K)+1 \}.
    \end{equation}
    Since $J$ (and $K$) is generated in degree $n-1$, the minimal monomial generators of $J\cap K$ are of degree at least $n-1$. Therefore by Lemma~\ref{lem:reg-max-degree}, for (\ref{reg-equa}) to hold, the ideal $J\cap K$ must be generated in degree $n-1$, just like $J$ and $K$ themselves. By Lemma~\ref{lem:formulas}, we must have that $J$ contains or is contained in $K$, as desired.

    \emph{(ii) $\implies $ (i):} Assume that $J$ and $K$ have linear resolution or linear quotients, and one contains the other. Without loss of generality, assume that $J\subseteq K$. For this direction, we will deal with the hypothesis of having linear resolution and linear quotients separately.
    
    Assume that $J$ and $K$ have linear resolution. Thus by Lemma~\ref{lem:x-splitting-special}, $I=x_nJ+y_nK$ is an $x_n$-splitting. Therefore, we have
    \begin{align*}
        \reg I &= \max\{ \reg J+1 ,\reg K +1, \reg (J\cap K) +1 \} \\
        &= \max\{ \reg J+1 ,\reg K +1, \reg (J) +1 \} \\
        &=n,
    \end{align*}
    i.e., $I$ has linear resolution, as desired.

    Now assume that $J$ and $K$ have linear quotients. Since $J\subseteq K$ and $J$ and $K$ are generated in the same degree, we have $\mingens(J)\subseteq \mingens(K)$. Since $K$ has linear quotients, we can set $\mingens(K)=\{m_1,\dots, m_p,m_{p+1},\dots, m_q\}$ such that
    \begin{equation}\label{LQ-J-K}
        (m_1,\dots, m_i)\colon m_{i+1} \text{ is generated by a subset of variables for any } i\in [q-1],
    \end{equation}
    and that $\mingens(J)=\{m_1,\dots, m_p\}$. Note that we now have
    \[
    \mingens(I)=\mingens(x_nJ)\sqcup \mingens(y_nK) = \{ x_nm_i, y_nm_j \mid i\in [p] \text{ and }j\in [q] \}.
    \]
    It now suffices to show that
    \begin{equation}\label{LQ-1}
        (y_nm_1,\cdots, y_nm_j)\colon y_nm_{j+1} \text{ is generated by variables for any } j\in [q-1]
    \end{equation}
    and 
    \begin{equation}\label{LQ-2}
        (y_nm_1,\dots, y_nm_q, x_nm_1,\dots,  x_nm_{i-1} )\colon x_nm_{i} \text{ is generated by variables for any } i\in [p].
    \end{equation}
    For the colon formulas below, we refer to Lemma~\ref{lem:formulas}. Indeed, (\ref{LQ-1}) follows immediately from (\ref{LQ-J-K}) as $(y_nm_1,\cdots, y_nm_j)\colon y_nm_{j+1}=(m_1,\cdots, m_j)\colon m_{j+1}$ for any $j\in [q-1]$. As for (\ref{LQ-2}), first note that
    \[
    (y_nm_1,\dots, y_nm_q)\colon x_nm_{1} =(y_n).
    \]
    Next, an application of Lemma~\ref{lem:formulas} gives
    \begin{align*}
        &\ \ \ (y_nm_1,\dots, y_nm_q, x_nm_1,\dots,  x_nm_{i-1} )\colon x_nm_{i} \\
        &= (y_nm_1,\dots, y_nm_q)\colon x_nm_{i} +(x_nm_1,\dots,  x_nm_{i-1} )\colon x_nm_{i}  \\
        &= (y_n) + (\text{a subset of } \{x_1,\dots, x_{n-1},y_1,\dots, y_{n-1}\})
    \end{align*}
    for any $i\in [2,p]$. The statement (\ref{LQ-2}) then follows. This concludes the proof.
\end{proof}

The preceding theorem is essentially a recursive formula to verify whether a polarized neural ideal on $n$ neurons  generated in degree $n$ has linear resolution or linear quotients. Basically, the two properties follow the same process of verification. The next result, from this perspective, is hardly surprising.

\begin{theorem}\label{thm:LR-LQ-the-same}
    Let $I$ be a polarized neural ideal on $n$ neurons and assume that $I$ is generated in degree $n$. Then $I$ has linear resolution if and only if it has linear quotients.
\end{theorem}

\begin{proof}
    We can set $I=x_nJ+y_nK$ where $J,K$ are (uniquely determined) ideals generated by monomials that do not involve $x_n$ or $y_n$. In particular, $J,K$ are neural ideals on $n-1$ neurons, and generated in degree $n-1$.  By Theorem~\ref{thm:linear-resolution/quotients-inductive}, $I$ has linear resolution or linear quotients if and only if so are $J$ and $K$, and either $J\subseteq K$ or $K\subseteq J$. In other words, whether the properties of having linear resolution and having linear quotients are equivalent for the neural ideal $I$ on $n$ neurons entirely depends on whether they are equivalent for the neural ideals $J$ and $K$ on $n-1$ neurons. The result then follows from induction and the fact that these two properties are equivalent for any monomial ideal generated in degree 1 (a subset of variables).   
\end{proof}

One may ask, how about only assuming that $I$ is generated in degree $d$ for some $d\leq n$? For special values of $d$ and $n$, we are confident that an analog can be obtained, but not for all pairs of $(d,n)$. The reason is that \emph{any} squarefree monomial ideal in $n$ variables $x_1,\dots, x_n$ is a polarized neural ideal on $n$ neurons, with the variables $y_1,\dots, y_n$ simply never appearing in the generators. Classifying polarized neural ideals with linear resolution (or linear quotients), is therefore tantamount to doing so for all squarefree monomial ideals, a task that has been studied intensively. It is also worth noting that from this perspective, having linear resolution and having linear quotients are not equivalent in general for polarized neural ideals (see \cite{LR-not-LQ}).

\bibliographystyle{amsplain}
\bibliography{refs}
\end{document}